\documentclass[12pt]{article}
\usepackage{latexsym}
\usepackage{amssymb}
\author{Mark Ramras\\ Department of Mathematics\\
Northeastern University\\ Boston, MA 02115}
\title{Routing Linear Permutations on Fibonacci and Lucas Cubes}
\begin{document}
\maketitle
\begin{abstract}
\indent\indent  In recent years there has been much interest in certain subcubes of hypercubes, namely Fibonacci cubes and Lucas cubes (and their generalized versions).  In this article we consider off-line routing of linear permutations on these cubes.  The model of routing we use regards edges as bi-directional, and we do not allow queues of length greater than one.  Messages start out at different vertices, and in movements synchronized with a clock, move to an adjacent vertex or remain where they are, so that at the next stage there is still exactly one message per vertex.  This is the routing model we defined in an earlier paper.
\end{abstract}

\begin{section}{Introduction}
  It is well-known that the $n$-dimensional hypercube $Q_n$ is the graph whose vertices are all binary strings of length $n$, and whose edges join those pairs of vertices which differ in exactly one position.  The Fibonacci cube of dimension $n$, introduced in \cite{Hsu}, which we shall denote by $F_n$ is the subgraph consisting of those strings with no two adjacent ones.  It is easy to see that $|\,V(F_n)\,|=f_n$, the $n^{\rm th}$ Fibonacci number.  For a comprehensive survey of the structure of Fibonacci cubes, see \cite{Klav}.
  The Lucas cube of dimension $n$, which we shall denote by $L_n$ is defined analogously, except that now the first and last positions, $x_1$ and $x_n$ of the string $(x_1,x_2,\ldots,x_n)$ are considered adjacent, and so may not both be $1$.  $|\,V(L_n)\,|=l_n,$ the $n^{\rm th}$ Lucas number.
  
In \cite{KlavMoll}  the automorphism group of $F_n$ is computed (it is $\mathbb{Z}_2$) and in \cite{Salv} that of $L_n$ (it is the dihedral group $D_{2n}$).  In the next section we determine the linear permutations of $F_n$ and $L_n$, {\it i.e.} the $n\times n$ $0-1$ matrices $A$ such that for all $x\in F_n ({\rm resp.} L_n), Ax\in F_n ({\rm resp.} L_n)$, and identify the groups they form.  In section 3 we give routings for these permutations, presented as products of permutations of $F_n ({\rm resp.} L_n)$, that move messages to adjacent vertices or leave them fixed, following the protocol and method of \cite{Ram}.

\end{section}

\begin{section}{Identifying the linear permutations of $F_n$ and of $L_n$}

\begin{subsection}{The linear permutations of $F_n$}
\newtheorem{definition}{Definition}  Let $A$ be an $n\times n$ invertible binary matrix.  We say that $A$ is $F_n$-good (resp. $L_n$-good) if $A(F_n)\subseteq F_n$ (resp. $A(L_n)\subseteq L_n$).

We first determine the $F_n$-good matrices.  From now on, assume that $A$ is an $n\times n$ invertible binary matrix, and assume that $A_n$ is $F_n$-good.  $A_{i,j}$ denotes the entry in row $i$, column $j$, $A_i$ the $i^{\,{\rm th}}$ row, and $A^{(j)}$ the $j^{\,{\rm th}}$ column.

\newtheorem{lemma}{Lemma}

\begin{lemma}
Suppose $A$ is $F_n$-good and has a $1$ in row $i$, column $j$.  Then a $1$ in row $i-1$ or row $i+1$ can only occur in columns $j-1$ or $j+1$.
\end{lemma}

{\em Proof.}  Suppose row $i+1$ has a $1$ in column $k$.  By definition of $F_n$, if $\vec{e}_j$ denotes the $j^{\rm th}$ standard basis vector of $\mathbb{Z}_2^n$, then $A\vec{e}_j=$ column $j$ does not have adjacent $1$'s and so $A_{i-1,j}=0=A_{i+1,j}$.  Suppose $k<j-1$ or $k>j+1$ and $A_{i-1,k}$ or $A_{i+1,k}$ or $A_{i+1,k}=1$.  If either $A_{i-1,k}=1$ or $A_{i+1,k}=1$ then $A_{i,k}=0$.  In either case,
 $A(\vec{e}_j+\vec{e}_k)=A\vec{e}_j+A\vec{e}_k=A^{(j)}+A^{(k)}$ has two adjacent ones, contradicting the fact that $A$ is $F_n$-good unless $|\;k-j\;|=1,$ in which case $\vec{e}_j+\vec{e}_k\notin F_n$.  So $k=j-1$ or $j+1$.  \hfill   $\Box$

\begin{lemma}
 Suppose $A$ is $F_n$-good, where $n\geq 3$.
If row $i$ of $A$ has two $1$'s, occurring in columns $j$ and $k$, with $j<k$, then $k=j+2$ and $i=1$ or $n$.
\end{lemma}

{\em Proof.}  By Lemma 1, because of the $1$ in row $i$, column $j$, a $1$ in row $i-1$ (if $i\geq 2$) must occur occur in either column $j-1$ or column $j+1$, and because of the $1$ in row $i$, column $k$, any $1$ in row $i-1$ must also occur in either column $k-1$ or column $k+1$.  To satisfy both of these conditions, since $j<k$ we must have $j+1=k-1$, {\it i.e.} $k=j+2$.  If $i\leq n-1$, the same argument holds for row $i+1$.  So if $2\leq i\leq n-1$, row $A_{i-1}=\vec{e}_{j+1}=$ row $A_{i+1}$, contradicting the assumption that $A$ is invertible.  Hence $i=1$ or $n$.  \hfill   $\Box$

\newtheorem{corollary}{Corollary}
\begin{corollary}
Suppose $A$ is $F_n$-good.  Then\\
$(1)$  Any $1$ in column $1$ must occur in row $1$ or row $n$.  Similarly, a $1$ in column $n$ must occur in row $1$ or row $n$.\\
$(2)$  $A_2, A_3, \ldots, A_{n-1}$ is a permutation of $n-2$ of the rows of the identity matrix $I_n$.\\
$(3)$  Each of $A_1$ and $A_n$ has either one or two $1$'s.\\
$(4)$  In row $1$ there must be a $1$ in either column $1$ or column $n$.  The same is true for row $n$.  A second $1$ $($if there is one$)$ in rows $A_1$ or $A_n$  must occur in column $A^{(3)}$ or $A^{(n-2)}$.
\end{corollary}

{\em Proof.}  (1)  Suppose $A_{i,1}=1$, and $2\leq i\leq n-1.$  Then since $A^{1}\in F_n$, $A_{i-1,1}=A_{i+1,1}=0.$  By Lemma 1, the only 1 in $A_{i-1}$ and the only 1 in $A_{i+1}$ both occur in column $A^{(2)}$.
Hence $A_{i-1}=\vec{e}_2=A_{i+1},$ contradicting the invertibility of $A$.  Thus for $2\leq i\leq n-1, A_{i,1}=0$.  Since $A^{(1)}\neq \vec{0}$, either $A_{1,1}=1$ or $A_{n,1}=1.$  The argument for column $A^{(n)}$ is analogous.
\vspace{.1in}

(2)  By Lemma 2, for $2\leq i\leq n-1$, $A_i$ has exactly one $1$, and two rows $A_{i_1}$ and $A_{i_2}$ are equal if their $1$'s occur in the same column.  Since $A$ is invertible, this cannot happen.  Thus for $2\leq i\leq n-1$, $A_2, A_3, \ldots, A_{n-1}$ is a permutation of $n-2$ of the rows of the identity matrix $I_n$.
\vspace{.1in}

(3)  This is part of Lemma 2, since for integers $j, k,$ and $l$ with $j<k<l$ not all 3 of their differences can be 2.
\vspace{.1in}

(4)  By (1), either $A_{1,1}=1$ or $A_{1,n}=1$.  By Lemma 2, if both are 1, then $n=1+2=3$.  If $n\geq 4$, then either $A_1=\vec{e}_1$, or $A_1=\vec{e}_1+\vec{e}_3$, or $A_1=\vec{e}_n$ or $A_1=\vec{e}_{n-2}+\vec{e}_n$.  The same statements hold for $A_n$, for analogous reasons, except that $A_1\neq A_n$.  \hfill  $\Box$

\begin{corollary}
If $A_{1,1}=1$, then for $2\leq i\leq n-1, A_i=\vec{e}_i$.\\
If $A_{1,n}=1,$ then for $2\leq i\leq n-1, A_i=\vec{e}_{n-i+1}$.
\end{corollary}

{\em Proof.}  Suppose $A_{1,1}=1$.  Then $A_{2,1}=0$.  By Lemma 1, $A_2=\vec{e}_2$.  Let $2\leq i\leq n-2$ and assume, inductively, that $A_i=\vec{e}_i$.  Then again by Lemma 1, $A_{i+1}=\vec{e}_{i+1}$.  So by induction, for $2\leq i\leq n-1, A_i=\vec{e}_i$.  \\
The argument when $A_{1,n}=1$ is entirely analogous.  \hfill   $\Box$

\vspace{.2in}

Next, we shall exhibit the group of $F_n$-good linear permutations.
\newtheorem{theorem}{Theorem}
\begin{theorem}
For $n\neq 3,$ there are precisely $8$ $F_n$-good linear permutations and they form the dihedral group $D_4$.  For $n=3, I+E_{1,3}+E_{3,1}=C+E_{3,3}+E_{1,1}$ is not invertible $($rows $1$ and $3$ are equal$)$, and so $F_3$ has exactly $6 F_3$-good linear permutations, and they form the permutation group $\frak{S}_3$.
\end{theorem}

{\em Proof.}  Denote by $E_{i,j}$ the $n\times n$ matrix whose single non-zero entry is a $1$ in the $(i,j)^{\,{\rm th}}$ position.  Let $C$ denote the matrix such that for $1\leq i\leq n, C_{i,n-i}=1$ and for $1\leq j\leq n,
C_{i,j}=0$ if $i+j\neq n$.  Suppose that $n\neq 3$.  We claim that the following set of 8 matrices constitutes the set of $F_n$-good matrices:  
$$\{I, I+E_{1,3}, I+E_{n,n-2}, I+E_{1,3}+E_{n,n-2}\}$$
$$\bigcup \{C, C+E_{n,3}, C+E_{1,n-2}, C+E_{n,3}+E_{1,n-2}\}.$$
From the previous corollary we see that any $F_n$-good matrix must be one of these $8$.  Next, we show that each {\em is} $F_n$-good.  Let $\vec{x}=(x_1,x_2,\ldots,x_n)$.  $C\vec{x}=(x_n,x_{n-1},\ldots,x_2,x_1)$.
If $C\vec{x}\notin F_n$ then for some $i$ with $n\geq i\geq 2, x_i$ and $x_{i-1}$ are both $1$.  But then $\vec{x}\notin F_n$.  Thus $C$ is $F_n$-good.  Next, consider $A=I+E_{1,3}$.  $A\vec{x}=(x_1+x_3,x_2,x_3,x_4,\ldots,x_n)$.  Since $\vec{x}\in F_n$, if $A\vec{x}\notin F_n$ then $x_1+x_3=1$ and either $x_1=1$ or $x_2=1.$  For $x_1+x_3$ to be 1, exactly one of $x_1$ and $x_3$ is 1.  In either case, among the first three components of $\vec{x}$, two adjacent ones are 1.  Thus $\vec{x}\notin F_n$.  This contradiction shows that $I+E_{1,3}$ is $F_n$-good.  We'll demonstrate the proof for one more example:
  $A=C+E_{1,n-2}+E_{n,3}$.  $A\vec{x}=C\vec{x}+E_{1,n-2}\vec{x}+E_{n,3}\vec{x}=(x_n+x_{n-2},x_{n-1},x_{n-2},\ldots,x_3,x_2,x_1+x_3)$.  If $A\vec{x}\notin F_n$ then either $x_n+x_{n-2}=x_{n-1}=1$ or $x_2=x_1+x_3=1$.  In the first case, exactly one of $x_n$ and $x_{n-2}$ is 1, and so together with $x_{n-1}$ we have two adjacent 1's in $\vec{x}$, or, in the second case, exactly one of $x_1$ and $x_3$ is 1, and so together with $x_2$ we again have two adjacent 1's in $\vec{x}$.  In either case, $\vec{x}\notin F_n$.  Hence $A$ must be $F_n$-good.
  
  Next, we must show that these $8 F_n$-good matrices form the dihedral group $D_4$.  It is an easy (but tedious) check to see that this set is closed under multiplication.  Since it is also composed of invertible matrices, it forms a group.  To see that this group is the dihedral group $D_4$ we will exhibit matrices $A$ and $B$ such that $A^4=B^2=I$, $A^2\neq I$, and $BAB=A^3$.  Let $A=C+E_{1,n-2}$ and $B=I+E_{1,3}$.
Since $C^2=I$ and $n\neq 3$, $A^2=I+CE_{1,n-2}+E_{1,n-2}C$.  Then $A^2=I+E_{1,3}+E_{1,n-2}$, and so $A^4=I$.  Clearly $B^2=I$.  Now $A^{-1}=A^3=(C+E_{1,n-2})(I+E_{1,3}+E_{1,n-2})=C+CE_{1,n-2}+E_{1,3}+E_{1,n-2}=C+CE_{1,n-2}=C+E_{1,n-2}$.  On the other hand, $BAB=(I+E_{1,3})(C+E_{1,n-2})(I+E_{1,3})=(C+CE_{1,3}+E_{1,n-2})(I+E_{1,3})=(C+E_{1,n-2}+E_{1,n-2})(I+E_{1,3})= C(I+E_{1,3}=C+CE_{1,3}=C+E_{1,2}.$  Hence $BAB=A^{-1}$.

For $n=3,$ the first and third rows of $I+E_{1,3}+E_{3,1}$ are equal and so the matrix is not invertible.  The same is true for $C+E_{n,3}+ E_{1,n-2}$.  Thus the order of $G_3$ is 6.  Finally, $G_3$ is non-abelian, since, for example, $(I+E_{1,3})(I+E_{3,1})=I+E_{1,3}+E_{3,1}+E_{11}\neq I+E_{3,1}+E_{1,3}+E_{3,3}=(I+E_{3,1}(I+E_{1,3}$.  Hence $G_3\cong \frak{S}_3$.  \hfill    $\Box$

\end{subsection}
\begin{subsection}{The linear permutations of $L_n$}

In $L_n$ entries in the first and last postions are considered adjacent, in addition to those in positions $i$ and $i+1$, for $1\leq i\leq n-1$.  Thus $\vec{x}\in L_n \Longleftrightarrow \vec{x}\in F_n$ and not both $x_1$ and $x_n$ are $1$.  Thus for an $n\times n$ binary matrix $A$ that is invertible, $A$ is $L_n$-good provided that for all $\vec{x}\in \mathbb{Z}_2^n$ with no $1$'s in adjacent positions, $A\vec{x}$ also has no $1$'s in adjacent positions.

We will determine the $L_n$-good matrices  $A$ as we did for the $F_n$-good matrices.\\

{\bf Remark}  Unlike $F_n$, $L_n$ has the property that if $\vec{x}\in L_n$ and $\vec{y}$ is obtained from $\vec{x}$ by a cyclic permutation of the coordinates of $\vec{x}$, then $\vec{y}\in L_n$. 

\begin{corollary}
If $A$ is an $L_n$-good matrix, then so is any matrix obtained from $A$ by a cyclic permutation of its rows.  In particular, $C$ is $L_n$-good.
\end{corollary}

{\em Proof.}  Let $B$ be the matrix such that for $1\leq i\leq n-1, B_i=A_{i+1}$, and $B_n=A_1$.  Then $B=A\left(\vec{e_2},\vec{e_3},\ldots ,\vec{e_n},\vec{e_1}\right)^T$, where $T$ denotes the transpose.  Let $\vec{x}\in L_n$.  So $B\vec{x}=(A\left(\vec{e_2},\vec{e_3},\ldots ,\vec{e_n},\vec{e_1}\right)^T)\vec{x}=A(\left(\vec{e_2},\vec{e_3},\ldots ,\vec{e_n},\vec{e_1}\right)^T)\vec{x}=A\vec{y},$ where $\vec{y}=(x_2,x_3,\ldots,x_n,x_1)^T$.  Since $\vec{x}\in L_n$, by the Remark, $\vec{y}\in L_n$.  But $A$ is $L_n$-good, so $A\vec{y}\in L_n$, {\it i.e.} $B\vec{x}\in L_n$.  Hence $B$ is $L_n$-good.  The second statement follows from the fact that $C$ is obtained from $I$ by a cyclic permutation of its rows.  \hfill  $\Box$

\begin{lemma}
If $A$ is $L_n$-good then no column of $A$ has $1$'s in both the first and last rows.
\end{lemma} 

{\em Proof.}  $\vec{e_1}+\vec{e_n}\notin L_n$ since the first and last positions are considered adjacent.  Thus for all $j$, $A\vec{e_j}\neq \vec{e_1}+\vec{e_n}$ and so $A^{(j)}$ does not have $1$'s in its first and last positions.  \hfill   $\Box$

\begin{lemma}
If $A$ is $L_n$-good and $A_{i,j}=1$ then a $1$ in row $A_{i-1}$ or $A_{i+1}$ can occur only in column $A^{(j-1)}$ or $A^{(j+1)}$.
\end{lemma}

{\em Proof.}  This is the same as Lemma 1 for $F_n$-good matrices, except that addition of subscripts is modulo $n$, so that when $i=1$, or $j=1$, by $i-1$ and $j-1$ we mean $n$, and if $i=n$ or $j=n$, by $i+1$ or $j+1$ we mean $1$.  So, for example, suppose that $A_{1,1}=1$.  We shall show that a 1 in $A_n$ occurs only in column $A^{(n)}$ or column $A^{(2)}$.  We know from Lemma 3 that $k\neq 1$.  Suppose that $A_{n,k}=1$ for some $k$ with $3\leq k\leq n-1$.  $\vec{x}=\vec{e}_1+\vec{e}_k$ has no 1's in adjacent positions, so it belongs to $L_n$.  But $A\vec{x}=A^{(1)}+A^{(k)}$ has 1 in both the first and last positions, so that $A\vec{x}\notin L_n$.  This contradiction means that if $A_{n,k}=1$ then $k=2$ or $n$.  The cases $A_{1,n}=1, A_{n,1}=1,$ and $A_{n,n}=1$ are handled similarly.  \hfill   $\Box$ 
\vspace{.2in}

Next, we have the analogue of Lemma 2, for $L_n$.
\begin{lemma}
Suppose that $A$ is $L_n$-good, and $A_{i,j}=1=A_{i,k}$, where $j<k$.  \\
$(i)$  If $2\leq j$ and $k\leq n=1,$ then $k=j+2$ and $A_{i+1,j+1}=1$.  Thus $A_{i+1}=\vec{e}_{j+1}$.\\
$(ii)$  If  $j=1$ then  $k=3$ and $A_{i+1,2}=1$, so that $A_{i+1}=\vec{e}_2$, or $k=n-1$ and $A_{i+1,n}=1$, and hence $A_{i+1}=\vec{e}_n$.\\
$(iii)$  If  $k=n,$ then $j=2$ and $A_{i+1,1}=1$, so that $A_{i+1}=\vec{e_1}$, or $j=n-2$ and $A_{i+1,n-1}=1$, and thus $A_{i+1}=\vec{e}_{n-1}.$\\
Similarly, in row $A_{i-1}$ a $1$ can occur only in column $j+1$.  Thus $A_{i-1}=\vec{e}_{j+1}$.  Hence every row of $A$ has exactly one $1$.   
\end{lemma}

{\em Proof.}  (i)  By Lemma 4, if $2\leq j\leq n-2$ and $j<k\leq n-1$, then since $A_{i,j}=1$, a $1$ in row $A_{i+1}$ can occur only in column $A^{(j-1)}$ or column $A^{(j+1)}$.  Similarly, since $A_{i,k}=1$, a $1$ in row $A_{i+1}$ can occur only in column $A^{(k-1)}$ or column $A^{(k+1)}$.  Thus $j+1=k-1$ and so $k=j+2$.  Furthermore, since $A$ is invertible, {\it some} entry in row $A_{i+1}$ must be $1$, we must have $A_{i+1,j+1}=1$, and so $A_{i+1}=\vec{e}_{j+1}$.  The same argument holds for row $A_{i-1}$.  Therefore $A_{i-1}=A_{i+1}$, contradicting the invertibility of $A$.  Thus every row of $A$ has exactly one $1$, and $A$ is a permutation of the identity matrix $I$. \\
The proofs of $(ii)$ and $(iii)$ are similar.  \hfill  $\Box$

\begin{theorem} $($Classification of $G_n=$ the group of $L_n$-good matrices$)$\\
$(i)$  $G_2$ is the group of order $2$, consisting of $I$ and $C=\left(\begin{array}{cc} 0&1\\ 1&0\end{array}\right)$.\\
$(ii)$ $G_3$ is the group of $3\times 3$ permutation matrices, {\em i.e.} the dihedral group $D_3$ $($or equivalently, the symmetric group ${\cal S}_3$$)$.\\
$(iii)$  $G_4$ is the group of $4\times 4$ matrices consisting of $I$, $I+E_{1,3}$, $I+E_{3,1}$,  $I+E_{2,4}$, $I+E_{1,3}+E_{2,4}$,  $I+E_{4,2}$, $I+E_{1,3}+E_{4,2}$, $C$, $C+E_{1,2}$, $C+E_{2,1}$, $C+E_{1,2}+E_{2,1}$, $C+E_{3,4}$, $C+E_{4,3}$, $C+E_{2,1}+E_{3,4}$,  $C+E_{1,2}+E_{3,4}$, $C+E_{1,2}+E_{4,3}$ and all matrices obtained from these by cyclic permutations of the rows.  The order of $G_4$ is $72$.\\
$(iv)$  For $n\geq 5$, each $L_n$-good matrix has exactly one $1$ in each row and exactly one $1$ in each column.  $G_n$ consists of all the cyclic permutations of the rows of $I$ and so $G_n\cong \mathbb{Z}_n$.
\end{theorem}

{\em Proof.}  $(i)$  None of the other $2\times 2$ matrices is invertible.\\
$(ii)$  For $n=3$, any two positions in $\vec{x}$ are adjacent.  If, say, the $j^{{\rm th}}$ column contains two $1$'s then $A\vec{e}_j\notin L_3$, and $A$ is not $L_3$-good.  Thus each column of $A$, and therefore each row of $A$, contains exactly one 1.  So $A$ is a permutation matrix, {\it i.e.} is obtained from $I$ by a permutation of its rows.   If $\pi$ is a non-cyclic permutation, then $\pi$ is a transposition.  But for $n=3$, any two rows may be considered adjacent.  So with no loss of generality, we may assume that $\pi=(1,2).$  Then $\pi (I)\vec{x}=(x_2,x_1,x_3)$.  Let $\vec{x}\in L_3$.  If any two 1's occur
in $(x_2,x_3,x_1)$ then $\vec{x}$ has two adjacent 1's, which is a contradiction.  Thus $\pi (I)\vec{x}\in L_3$ and so $\pi(I)$ is $L_3$-good.  Since the symmetric group is generated by the transpositions $G_3\cong \frak{S}_3$.\\
$(iii)$  It is easy to check that each of these matrices is $F_4$-good.  By the {\it weight} of a row we mean the number of 1's in it.  There are three cases for an $F_4$-good matrix $A$:  (1)  no row of $A$ has weight 2, (2) exactly one row has weight two, or (3)  exactly two rows have weight two.  In case (1) we have the 8 cyclic permutations of $I$.  In cases (2) and (3), it follows from Lemma 5, ($i$) and ($ii$) that there are just two choices for the row of weight 2:  $[1010]$ or $[0101]$.  Now in case (2), there are 4 possible positions for the row of weight 2.  Since the assume that  $[1010]$ is row 1 of $A$.  Then $A_{2,1}=A_{2,3}=A_{4,1}=A_{4,3}=0$.  $A_2$ must be either $\vec{e_2}$ or $\vec{e_4}$, and $A_4$ must be either $\vec{e_4}$ or $\vec{e_2}$ (and $A_2\neq A_4$).  Thus there are 4 $L_4$-good matrices with $A_1=[1010]$.  Since $L_n$-good matrices are closed under cyclic permutations of rows, we have $4\times 4$ $L_4$-good matrices with the single row of weight 2 being [1010].  Similarly, there are $4\times 4$ $L_4$-good matrices whose single row of weight 2 is [0101].  Thus, in case (2) we have $2\times 4\times 4 = 32$ $L_4$-good matrices.

Case (3) has 2 subcases:  (a)  the 2 rows of weight 2 are consecutive, and (b) they are not.  \\
For subcase (a):  first suppose that $A_1=[1010]$ and $A_2=[0101]$.  We obtain the following 4 $L_4$-good matrices:
\vspace{.4in}

$\left(\begin{array}{llll}1 & 0& 1 & 0\\
                                  0 & 1 & 0 & 1\\
                                  1 & 0 & 0 & 0\\
                                  0 & 1 & 0 & 0
                                  \end{array}\right)$, $\left(\begin{array}{llll}1 & 0& 1 & 0\\
                                  0 & 1 & 0 & 1\\
                                  1 & 0 & 0 & 0\\
                                  0 & 0 & 0 & 1
                                  \end{array}\right)$, $\left(\begin{array}{llll}1 & 0& 1 & 0\\
                                  0 & 1 & 0 & 1\\
                                  0 & 0 & 1 & 0\\
                                  0 & 1 & 0 & 0
                                  \end{array}\right)$, $\left(\begin{array}{llll}1 & 0& 1 & 0\\
                                  0 & 1 & 0 & 1\\
                                  0 & 0 & 1 & 0\\
                                  0 & 0 & 0 & 1
                                  \end{array}\right)$.
  \vspace{.1in}

 Similarly, if we interchange the first two rows, we get the following 4 $L_4$-good matrices:
 \vspace{.4in}
 
 $\left(\begin{array}{llll}0 & 1& 0 & 1\\
                                  1 & 0 & 1 & 0\\
                                  0 & 1 & 0 & 0\\
                                  1 & 0 & 0 & 0
                                  \end{array}\right)$,   $\left(\begin{array}{llll}0 & 1& 0 & 1\\
                                  1 & 0 & 1 & 0\\
                                  0 & 1 & 0 & 0\\
                                  0 & 0 & 1 & 0
                                  \end{array}\right)$,     $\left(\begin{array}{llll}0 & 1& 0 & 1\\
                                  1 & 0 & 1 & 0\\
                                  0 & 0 & 0 & 1\\
                                  1 & 0 & 0 & 0
                                  \end{array}\right)$,    $\left(\begin{array}{llll}0 & 1& 0 & 1\\
                                  1 & 0 & 1 & 0\\
                                  0 & 0 & 0 & 1\\
                                  0 & 0 & 1 & 0
                                  \end{array}\right)$.  
      \vspace{.1in}
      
  Since there are 4 pairs of consecutive rows, there are a total of $(4+4)\times 4=32$ $L_4$-good matrices in case 3(a).

 For subcase (b):  If the 2 rows of weight 1 are not consecutive, then each of the rows between them must be the zero row, since no column can have two consecutive 1's.  Thus there are no $L_4$-good matrices of this type.
 
 Hence the total number of $L_4$-good matrices is $8+32+ 32+0=72$.

$(iv)$  Let $n\geq 5$, and suppose that for some $i$, $1\leq i\leq n$, row $A_i$ has at least two 1's.  Say that $A_{i,j}=1=A_{i,k}$, where $j<k$.  Then $A_{i+1}=\vec{e}_{j+1}=A_{i-1}$.  This contradicts the fact that $A$ is invertible.  So each row of $A$ is $\vec{e}_j$, for a unique $j$.  Thus $A$ is obtained from $I$ by a permutation of its rows.  We must show that this permutation is {\em cyclic}.  Suppose that $A_1=\vec{e}_j$.  Then $A_2=$ either $\vec{e}_{j+1}$ or $\vec{e}_{j+3}$.  First suppose it is $\vec{e}_{j+1}$.  We claim that for all $1\leq i\leq n, A_i=\vec{e}_{j+i-1}$ (remember that subscripts are computed mod $n$).  Assume, inductively, that for $1\leq q\leq i-1, A_q=\vec{e}_{j+q+1}$.  Then $A_i=$ either $\vec{e}_{j+i+1}$ or $\vec{e}_{j+i-1}$.  If it is $\vec{e}_{j+i-1}$, then $A_{i-2}=\vec{e}_{j+(i-2)+1)}=\vec{e}_{j+i-1}=A_i$, contradicting the fact that no two rows of $A$ can be equal.  Hence $A_i=\vec{e}_{j+i+1},$ and so, by induction, for all $1\leq q\leq n, A_q=\vec{e}_{j+q+1}$.  Thus $A$ is obtained from $I$ via the cyclic permutation $i\mapsto j+i-1$.  A similar inductive argument shows that if $A_2=\vec{e}_{j+3}$ then for all $1\leq q\leq n, A_q=\vec{e}_{j+q+1}$, and so $A$ is obtained from $I$ via the cyclic permutation $i\mapsto j+i+1$.  \hfill  $\Box$

\end{subsection}
\end{section}
\begin{section}{Routings}
For a connected graph $G$, if $\pi$ is a permutation of $G$, we define  $$t(\pi)={\rm max}\{d_G(\pi(x),x)\,|\,x\in G\}.$$  We then consider the group Perm$(G)$ to be the Cayley graph whose generating set is $\Delta =\{\pi\,|\,t(\pi)\leq1\}$.  A $t$-fold product of elements of $\Delta$ equal to the permutation $\sigma$ is said to be a ``$t$-step routing of $\sigma$.
As discussed in \cite{Ram}, since we consider each edge to be doubled, {\it i.e.} one in each direction, every element $\tau$ of Perm$(G)$ {\em is} a finite product of elements of this generating set, and such a factoring we call a {\em routing} of the permutation $\tau$.
We will be interested in the two cases, $G=F_n$ and $G=L_n$.  For an $n\times n$ matrix $A$ which is $G$-good, we define the permutation $\tau_A$ by $\tau_A(\vec{x})=A\vec{x}$.
\begin{subsection}{Routings of permutations of $F_n$}
\begin{lemma}
For $A=I, I+E_{1,3},$ and  $I+E_{n,n-2}, t(\pi_A)=1$.  For $n\neq 3, t(I+E_{1,3}+E_{n,n-2})=1$.  For $n=3, I+E_{1,3}+E_{n,n-2}=I+E_{1,3}+E_{3,1}$ is not $F_3$-good, nor is $C+E_{n,3}+ E_{1,n-2}$.
\end{lemma}

{\em Proof}.  $d(A\vec{x},x)={\rm weight}(A\vec{x}+x)$.  For $A=I, I\vec{x}+\vec{x}=\vec{0}$ and $\rm weight(\vec{0}=0$.  For $A=I+E_{1,3}, A\vec{x}=\vec{x}+x_3\vec{e_1}$.  Thus $A\vec{x}+x=x_3\vec{e_1}$, whose weight is $0$ if $x_3=0$ and $1$ if $x_3=1.$  Therefore $t(A)=1$.  If $A=I+E_{n,n-2}$ then $A\vec{x}+x=E_{n,n-2}\vec{x}=x_{n-2}\vec{e}_n$, so again, $t(A)=1$.  If $A=I+E_{1,3}+E_{n,n-2}$, then $A\vec{x}+\vec{x}=x_3\vec{e}_1+x_{n-2}\vec{e}_n$, whose weight is $2$.   For $n\neq 3, E_{1,3}E_{n,n-2}=0$, and so $A=(I+E_{1,3})(I+E_{n,n-2})$ is a 2-step routing of $A$.    \hfill   $\Box$ 

\begin{lemma}
If $\vec{x}, \vec{y}, \vec{z}$ is a path in $F_n$, then $(\vec{x},\vec{z})=(\vec{x},\vec{y})(\vec{y},\vec{z})(\vec{x},\vec{y})$ is a $3$-step routing of $(\vec{x},\vec{z})$.
\end{lemma}

\begin{corollary}
$C$ is obtained from $I$ by the row permutation \linebreak $(1,n)(2,n-1)\dots (k,k+1)$ if $n=2k$ and by $(1,n)(2,n-1)\dots (k,k+2)$ if $n=2k+1$.  Each row transposition $(i,j)$ $({\rm where} \,|j-i\,|>1)$ can be routed in $3$ steps.  
Hence $C$ can be routed in $3\lfloor{n/2}\rfloor$ steps.
\end{corollary}

\begin{corollary}
$C+E_{n,3}=(I+E_{n,n-2})C$ and $C+E_{n,n-2}=(I+E_{n,3})C$.  Hence each of these can be routed in $1+3\lfloor{n/2}\rfloor$ steps.
\end{corollary}
\end{subsection}
\begin{subsection}{Routings of permutations of $L_n$}
\begin{lemma}
Let $A$ be the matrix obtained from $I$ by the row permutation $(1,2,3,\ldots,n)$.  Then the permutation $\tau_A$ corresponds to the product of transpositions $(1,n)(1,n-1)\dots(1,3)(1,2)$.  Each transposition $(1,i)$
corresponds to and element of $\Delta$, and hence $\tau_A$ has an $n$-step routing.
\end{lemma}

\begin{corollary}
If $A$ is obtained from $I$ by {\it any} cyclic permutation of the rows of $I$, then $\tau_A$ can be routed in at most $n$ steps.
\end{corollary}

{\em Proof.}  Any cyclic permutation corresponds to a power of $(1,2,\ldots,n)$.  This, in turn, is a product of disjoint cycles.  Each cycle of length $k$ is the product of $k$ transpositions, and therefore can be routed in $k$ steps.  Since the cycles are disjoint, the transpositions in one cycle are disjoint from those in the other cycles.  Hence the routings of these cycles can be carried out simultaneously, and so the number of steps in the routing is the maximum length of a cycle in this product, and thus is at most $n$.  \hfill    $\Box$

\begin{lemma}
$C$ can be routed in $3\lfloor{n/2}\rfloor$ steps.
\end{lemma}

{\em Proof.}  The routing for $\tau_C$ as a permutation of $F_n$ given in Corollary 4 works equally well for $t_C$ as a permutation of $L_n$.

\begin{corollary}
If $A$ is obtained from $C$ by a cyclic permutation of the rows of $C$ then $\tau_C$ can be routed in at most $5n/2$ steps.
\end{corollary}

{\em Proof.}  Let $\pi$ be a cyclic permutation of the rows of $C$.  Then $A=\pi(C)=C\pi(I)$.  By Lemma 8 $\tau_{\pi(I)}$ has can be routed in at most $n$ steps, and then by Lemma 9 $\tau_C$ can be routed in an additional 
$3\lfloor{n/2}\rfloor$ steps, for a total of $n+3\lfloor{n/2}\rfloor\leq 5n/2$ steps.   \hfill   $\Box$
\end{subsection}
\end{section}

\end{document}